\documentclass{article}

\usepackage{amsmath,amssymb,amscd,graphicx,psfrag,epsf,amsthm
}
\usepackage[all]{xy}

\newtheorem{thm}{Theorem}[section]

\newtheorem{cor}[thm]{Corollary}
\newtheorem{prop}[thm]{Proposition}
\newtheorem{conj}[thm]{Conjecture}
\newtheorem{problem}[thm]{Problem}

\theoremstyle{definition}
\newtheorem{defn}[thm]{Definition}
\newtheorem{exmp}[thm]{Example}

\newtheorem{condition}[thm]{Conditions}
\newtheorem{construction}[thm]{Construction}
\newtheorem{properties}[thm]{Properties}

\theoremstyle{remark}
\newtheorem{rem}[thm]{Remark}

\numberwithin{equation}{section}

\newcommand{\cS}{{\mathcal {S}}}
\newcommand{\cU}{{\mathcal {U}}}

\newcommand{\cR}{{\mathcal {R}}}

\renewcommand{\k}[0]{{\bf k}}  

\renewcommand{\c}[0]{{\mathbb C}}  

\renewcommand{\o}[0]{{\mathcal O}}

\newcommand{\m}[0]{{\mathfrak m}}

\newcommand{\p}[0]{{\mathbb P}}

\newcommand{\map}[0]{\dasharrow}

\newcommand{\spec}[0]{\operatorname{Spec}}

\newcommand{\coker}[0]{\operatorname{coker}}    
    
\newcommand{\Hom}[0]{\operatorname{Hom}}

\newcommand{\sing}[0]{\operatorname{Sing}}

\newcommand{\hilb}[0]{\operatorname{Hilb}}

\newcommand{\br}[0]{{br}}
\newcommand{\Br}[0]{{Br}}
\newcommand{\sym}[0]{\operatorname{Sym}}
\newcommand{\defm}[0]{\operatorname{Def}}
\newcommand{\ev}[0]{{ev}}

\begin{document}

\title{Rationally connected varieties}

\author{Carolina Araujo 
  \thanks{This work was partially completed during the period the author was employed by the 
    Clay Mathematics Institute as a Liftoff Fellow} \\
  IMPA \\
  Rio de Janeiro, Brazil \\
  caraujo@impa.br
}


\date{}

\maketitle

\begin{abstract}
The aim of these notes is to provide an introduction to the theory of rationally connected varieties,
as well as to discuss a recent result by T. Graber, J. Harris and J. Starr.
\end{abstract}


\section{Introduction}\label{intro}

Consider the problem of classifying the ``simplest'' projective varieties. 
In dimension $1$ there is not much to say. 
Smooth projective curves are classified by their genus.
It is not hard to argue that $\p^1$ is the simplest projective curve.
(See \cite{kollar_simplest_alg_var} for a detailed discussion.)

In dimension $2$ classification was completed by Enriques at the beginning of the twentieth century.
Rational surfaces (i.e., surfaces birational to $\p^2$) form a distinguished class of surfaces.
This class is very well behaved from the point of view of classification.
To illustrate this, we state some nice properties of the class of rational surfaces.

\begin{properties}[Nice properties of rational surfaces]\label{properties_rat_surfaces} \

\begin{enumerate}
  \item (Deformation invariance.) 
    Let $\cS\to B$ be a connected family of smooth complex projective surfaces.
    If $\cS_{b_o}$ is rational for some $b_0\in B$, then $\cS_b$ is rational for every $b\in B$.
  \item (Numerical characterization.)
    Let $S$ be a smooth complex projective surface. 
    Then $S$ is rational if and only if $H^0(S,\omega_S^{\otimes 2})=H^1(S,\o_S)=0$.
  \item (Geometric criterion.)
    Del Pezzo surfaces (i.e., smooth complex projective surfaces $S$ for which $-K_S$ is ample) are rational.
  \item (Well behavior under fibration.)
    Let $S$ be a smooth complex projective surface. 
    Assume that there exists a morphism $S\to \p^1$ whose general fiber is rational.
    Then $S$ is rational. 
\end{enumerate}  
\end{properties}

Now let us move to higher dimensions. 
After looking at the $2$-dimensional case, one is naturally led to considering
the class of rational varieties (i.e., varieties birational to $\p^n$, $n\geq 1$).
The drawback is that
this class of varieties behaves very badly in many aspects.
Consider for instance Properties~\ref{properties_rat_surfaces}.
They all seem to fail in higher dimensions.

\medskip
\begin{enumerate}
  \item Deformation invariance is conjectured to fail in dimensions $\geq 3$.
  \item Let $X_3\subset \p^4$ be a smooth cubic hypersurface. Then
    $H^0(X_3,\omega_{X_3}^{\otimes m})=H^i(X_3,\o_{X_3})=0$ for every $m,i\geq 1$,
    but $X_3$ is not rational by \cite{clemens_griffiths}.
  \item Let $X_3$ be as above. Then $-K_{X_3}$ is ample but $X_3$ is not rational.
  \item There are examples of $\p^1$-bundles over $\p^2$ that are not rational.
\end{enumerate}
\medskip

These are some of the reasons why the class of rational varieties 
is not suitable for classification purposes.
As a substitute, the more general concept of \emph{rationally connected varieties} 
was introduced in \cite{campana} and \cite{kmm2}.
As we shall see in the next section, this bigger class of varieties is very well behaved.
It satisfies many nice properties, including the analogs of Properties~\ref{properties_rat_surfaces}
(except perhaps one direction of the numerical characterization, which is conjectured to hold as well). 
In this context, rationally connected varieties should be viewed as the right higher dimensional analogs
of rational curves and rational surfaces.

In section~\ref{section:RCV} we define rationally connected varieties and prove several 
important properties that they satisfy.
In section~\ref{section:MRC} we introduce the concept of maximal rationally connected fibration.
In section~\ref{section:GHS} we explain the proof of a recent result by Graber, Harris and Starr
that we use in sections~\ref{section:RCV} and \ref{section:MRC}.
In section~\ref{section:final} we state some further results and open problems.

We refer to \cite{kollar} for a fairly complete treatment of the theory of rational curves on varieties.
We also refer to \cite{debarre} for a nice introduction to the subject.


\section{Rationally connected varieties}\label{section:RCV}

There are many different ways of defining rationally connected varieties.
We refer to \cite{kmm2} for the proof that the conditions below are equivalent.

\begin{defn}\label{def:RCV}
Let $X$ be a smooth complex projective variety of positive dimension.
We say that $X$ is \emph{rationally connected} if the following equivalent conditions hold.
\begin{enumerate}
  \item[1.] Two general points of $X$ can be connected by a chain of rational curves.
  \item[2.] Any two points of $X$ can be connected by a rational curve.
  \item[3.] Any finite set of points in $X$ can be connected by a rational curve.
  \item[4.] There exists a morphism $f:\p^1\to X$ such that 
$$
f^*T_X\ \cong \ \bigoplus_{i=1}^{\dim X}\o_{\p^1}(a_i), \ \text{\ with all } \ a_i\geq 1. 
$$
\end{enumerate} 
\end{defn}

\begin{rem}\ 
\begin{enumerate}
  \item A point is considered to be rationally connected.
  \item By a general point of a variety $X$, we mean a point in some dense open subset of $X$.
    The notion of general point depends upon the choice of proper closed subset of $X$ to be avoided.
    This should be clear from the context. 
\end{enumerate}
\end{rem}

Condition 4 above says that rationally connectedness can be detected in the presence of a single 
rational curve. Curves satisfying this condition enjoy many nice properties and
are very important in the study of rationally
connected varieties.

\begin{defn}
Let $X$ be a smooth projective variety, and $C\subset X$ a rational curve.
Let $f:\p^1\to C$ be a surjective morphism.
We say that $C$ (or $f$) is a \emph{very free} curve if 
$$
f^*T_X\ \cong \ \bigoplus_{i=1}^{\dim X}\o_{\p^1}(a_i), \ \text{\ with all } \ a_i\geq 1. 
$$
\end{defn}

Notice that $C$ is very free if and only if $H^1(\p^1, f^*T_X(-2))=0$.
Hence, by the Semicontinuity Theorem (see \cite[III.12]{hartshorne}), a general 
deformation of a very free curve is still very free.

Condition 4 in Definition~\ref{def:RCV} seems at first very different in flavor 
from the previous ones. In order to understand the relation between them, 
one needs to know a little about deformation theory. 
Intuitively, if $C$ is a very free rational curve on $X$, then any two points of 
$X$ can be connected by a (possibly reducible) deformation of $C$ in $X$.
The next proposition makes this precise. We refer to \cite[II.3]{kollar} for a proof.

\begin{prop}\label{geometric_v.free}
Let $X$ be a smooth complex projective variety.
Let $H$ be an irreducible component of $\hilb(X)$ whose general member parametrizes 
a rational curve. 
Then the following conditions are equivalent.
\begin{enumerate}
  \item The general member of $H$ parametrizes a very free curve.
  \item Two general points of $X$ can be connected by a curve parametrized by $H$. 
\end{enumerate}
\end{prop}

\begin{rem}\label{remark_v_free}
We see from Proposition~\ref{geometric_v.free} that the following two conditions are equivalent to 
conditions 1--4 in Definition~\ref{def:RCV}.
\begin{enumerate}
   \item[5.] Two general points of $X$ can be connected by a very free rational curve.
     (In fact, we may even require that
     \emph{any} two points of $X$ can be connected by a very free rational curve.)
   \item[6.] For a general point $x\in X$ and a general tangent vector $\xi \in T_xX$,
     there exists a very free curve on $X$ passing through $x$ with tangent direction $\xi$.
\end{enumerate} 

If $\dim X\geq 3$, we can require that the very free curves in these conditions 
are smooth (see \cite[II.1.8]{kollar}).
This is not necessarily true if $\dim X=2$.
\end{rem}

Next we state a very useful feature of very free curves.
We refer to \cite[II.3.7]{kollar} for a proof.

\begin{prop}\label{avoiding_codim_2}
Let $X$ be a smooth projective variety and $S\subset X$ a subset of codimension at least $2$.
Let $C\subset X$ be a very free curve and $P\in C\setminus S$ any point.
Then a general deformation of $C$ passing through $P$ does not intersect $S$.
\end{prop}

Let us look at some easy examples.

\begin{exmp}
Set $X=\p^n$.
Any two points of $X$ can be connected by a line.
Thus $\p^n$ is rationally connected.
In fact, any line $l$ on $\p^n$ is very free.
Indeed, by  restricting the exact sequence 
$$
0 \ \to \ \o_X \ \to \ \o_X(1)^{\oplus n+1}\ \to \ T_X \ \to \ 0
$$
to $l$, we see that $T_X|_l\cong \o(2)\oplus \o(1)^{\oplus n-1}$.
\end{exmp}

\begin{exmp}
Let $X_3\subset \p^4$ be a smooth cubic hypersurface.
It is easy to see that $X_3$ contains a line through every point. 
Let $l\subset X_3$ be a line. 
Then $\deg T_{X_3}|_l=-K_{X_3}\cdot l=2$.
On the other hand, for a very free curve $f:\p^1\to X_3$ we must have $\deg f^*T_{X_3}\geq 4$.
Hence lines on $X_3$ are not very free.

Now let us consider conics in $X_3$.
Let $x_1$ and $x_2$ be general points in $X_3$. Let $l_{12}$ be the line on $\p^4$
joining $x_1$ and $x_2$. 
Since $x_1$ and $x_2$ are general points, 
$l_{12}$ is not contained in $X_3$. So it meets $X_3$ in a third 
general point $x_3\in X_3$. Let $l_3\subset X_3$ be a line through $x_3$.
Let $L$ be the $2$-plane spanned by $l_{12}$ and $l_3$.
Then $L\cap X_3$ is a cubic plane curve containing $l_3$ as an irreducible
component. So there is a conic $C$ such that  $L\cap X_3=C\cup l_3$.
Moreover, $x_1, x_2\in L\setminus l_3$, and thus  $x_1, x_2 \in C$. 
So $X_3\subset \p^4$ is rationally connected, and, 
by Proposition~\ref{geometric_v.free}, a general conic on $X_3$ is very free.
(Recall that $X_3$ is not rational!) 
\end{exmp}

Now let us look at some properties of rationally connected varieties.

\begin{properties}[Nice properties of rationally connected varieties]\label{properties_RCV} \ 

\begin{enumerate}
\setcounter{enumi}{-1}
  \item Rationally connectedness is a birational invariant.
  \item Rationally connectedness is invariant under smooth deformation. 
  \item If $X$ is rationally connected, then $H^0(X,(\Omega_X^1)^{\otimes m})=0$ for every $m\geq 1$.
  \item Fano varieties (i.e., smooth complex projective varieties $X$ for which $-K_X$ is ample) 
    are rationally connected. 
    In particular, smooth hypersurfaces of degree $d$ in $\p^n$ are rationally connected for $d\leq n$. 
\end{enumerate}  
\end{properties}

\begin{proof}[Sketch of proof]
Birational invariance follows from condition 1 of Definition~\ref{def:RCV}, 
as it only refers to a dense open subset of the variety.

The existence of a very free curve is an open condition in smooth families, while 
condition 1 of Definition~\ref{def:RCV} is a closed condition. 
Thus being rationally connected is deformation invariant.

Now let $X$ be a rationally connected variety and pick a global section 
$\alpha\in H^0(X,(\Omega_X^1)^{\otimes m})$, $m\geq 1$.
Let $f:\p^1\to X$ be a very free curve.
Since $f^*\Omega_X^1 \cong  \bigoplus_{i=1}^{\dim X}\o_{\p^1}(-a_i)$, with all $a_i\geq 1$,
we have that $H^0(\p^1,f^*(\Omega_X^1)^{\otimes m})=0$, and thus $f^*\alpha=0$.
In other words, the restriction of $\alpha$ to any very free curve is identically zero. 
This implies that $\alpha\equiv 0$, as $X$ is covered by very free curves.
We have proved that $H^0(X,(\Omega_X^1)^{\otimes m})=0$.

The fact that Fano varieties are rationally connected was established 
in \cite{campana} and \cite{kmm3}.
\end{proof}

While properties (0)--(2) above follow almost imediately from the definition,
and property (3) was proved shortly after it,
the analog of Property~\ref{properties_rat_surfaces}(4) remained open for almost one decade. 
It was proved recently by  Graber, Harris and Starr in \cite{GHS}, as a corollary of their 
main theorem:

\begin{thm}\cite{GHS} \label{ghs} 
Let $X$ be a smooth complex projective variety and $f:X\to B$ a surjective morphism onto a smooth curve. 
If the general fiber of $f$ is 
rationally connected, then $f$ has a section. \end{thm}

We explain the proof of this theorem in section~\ref{section:GHS}.

\begin{cor}\label{fibration}
    (Rationally connected varieties are well behaved under fibration.)
    Let $X$ be a smooth complex projective variety. 
    Assume that there exists a surjective morphism $f: X\to Y$ with 
    $Y$ and the general fiber of $f$ rationally connected.
    Then $X$ is rationally connected. 
\end{cor}

\begin{proof}
In order to prove that $X$ is rationally connected, we shall check condition 1
of Definition~\ref{def:RCV}.
Let $x_1, x_2\in X$ be general points.
Set $y_1=f(x_1)$, $y_2=f(x_2)$, $X_1=f^{-1}\{y_1\}$ and $X_2=f^{-1}\{y_2\}$. 
Since $x_1$ and $x_2$ are general, $X_1$ and $X_2$ are smooth and rationally connected. 
Since $Y$ is rationally connected, there exists a rational curve $C$ 
joining $y_1$ and $y_2$. 
Let $X_C$ be a desingularization of $f^{-1}C$. 
Then $f$ induces a map $f_C:X_C\to \p^1$ satisfying the hypothesis of Theorem~\ref{ghs}.
Thus $f_C$ has a section, which yields a rational curve $l\subset X$ meeting $X_1$ and $X_2$.
Let $z_1\in l\cap X_1$ and $z_2\in l\cap X_2$.
Since $X_1$ and $X_2$ are rationally connected,
there are rational curves $l_1\subset X_1$ and $l_2\subset X_2$ with 
$x_1, z_1\in l_1$ and $x_2, z_2\in l_2$.
So $l_1\cup l\cup l_2$ is a chain of rational curves connecting $x_1$ and $x_2$.
\end{proof}


\section{The maximal rationally connected fibration}\label{section:MRC}

Let $X$ be a smooth complex projective variety.
In this section we are interested in the problem
of measuring how much $X$ fails to be rationally connected.
In this context, consider the following equivalence relation on $X$:
$$
(x,y) \ \in \ \cR \ \iff \ \text{$x$ and $y$ can be connected by a chain
                                of rational curves.}
$$
One good solution to this problem would be to obtain a fibration $\varphi:X\to Y$ 
such that every fiber of $\varphi$ is an equivalence class of $\cR$.
Such a fibration, however, does not always exist. 
There are many examples of varieties, such as $K3$ surfaces,
containing countably many rational curves.
It is clear that there cannot be a fibration as above for such varieties.

It turns out that one can find a good substitute for this fibration 
if one ignores countably many proper closed subvarieties of $X$.
We make this precise in the next theorem, which was proved in 
\cite{campana} and \cite{kmm2}.

\begin{thm}\label{mrc_quotient}
Let $X$ be a smooth complex projective variety.
Then there exists a dense open subset $X^0$ of $X$, a normal variety $T^0$, and 
a proper surjective morphism $\varphi^0:X^0\to T^0$ such that 
\begin{enumerate}
   \item the general fiber of $\varphi^0$ is rationally connected, and 
   \item the \emph{very general} fiber
         of $\varphi^0$ is an equivalence class of $\cR$. 
\end{enumerate}
Moreover, this morphism is unique up to birational equivalence.
\end{thm}

\begin{rem} 
By a very general point of a variety $X$, we mean a point outside the union 
of countably many proper closed subvarieties of $X$.
\end{rem}

We call the morphism $\varphi^0$ (or rather the birational class of $\varphi^0$) the 
\emph{maximal rationally connected fibration of $X$}.
We call $T^0$ (or rather the birational class of $T^0$)
the \emph{maximal rationally connected quotient of $X$}.
The terminology is explained by the following property of $\varphi^0$.
Suppose $\psi:X'\to Z'$ is another proper morphism from 
a dense open subset of $X$ onto a normal variety 
satisfying condition 1 of Theorem~\ref{mrc_quotient}. 
Then there exists a rational map $\rho:Z'\map T^0$ such that
$\varphi^0=\rho \circ \psi$.

The maximal rationally connected quotient of $X$ is a point 
if and only if $X$ is rationally connected.
It is $X$ itself if and only if $X$ is not uniruled. 
(A complex variety $X$ is uniruled if it contains a rational curve through every point.)
In \cite{kmm2}, the authors raised the question of whether a 
maximal rationally connected quotient could possibly be uniruled. 
This was finally settled as a corollary of 
Graber, Harris and Starr's Theorem~\ref{ghs}.

\begin{cor}[of Theorem~\ref{ghs}]
Let $X$ be a smooth complex projective variety.
Then the maximal rationally connected quotient of $X$ is not uniruled.
\end{cor}

\begin{proof}
Let $\varphi^0:X^0\to T^0$ be the maximal rationally connected fibration of $X$.
Let $t\in T^0$ be a very general point.
If $T^0$ is uniruled, then $T^0$ contains a rational curve $C$ through $t$.
Consider the restriction of $\varphi^0$ to $(\varphi^0)^{-1}C$.
After compactification and desingularization, 
we can apply Theorem~\ref{ghs}. 
So we get a rational curve on $X$ meeting the fiber $(\varphi^0)^{-1}\{t\}$
but not contained in it. 
This contradicts the assumption that $t$ is a very general point of $T^0$,
and thus $(\varphi^0)^{-1}\{t\}$ is an equivalence class of $\cR$.
\end{proof}


\section{A theorem by Graber, Harris and Starr}\label{section:GHS}

The goal of this section is to explain the proof of Graber, Harris and Starr's Theorem~\ref{ghs}
as given in \cite{GHS}.
Here is the general idea of the proof.
Let $f:X\to B$ be a surjective morphism from a smooth complex projective variety 
onto a smooth curve. 
Assume that the general fiber of $f$ is rationally connected.  
Let  $C\subset X$ be a multisection of $f$. 
We would like 
the space of deformations of $C$ in $X$ to be big enough
so that we can degenerate $C$ into a reducible curve 
containing a section of $f$ as an irreducible component.
So we ask the following question.
{\it Starting with any multisection $C\subset X$, how can we make the space of deformations of
$C$ bigger?}
This is where the rationally connected fibers of $f$ come in. 
Let $P_1$ be a general point in $C$. The fiber of $f$ containing $P_1$ is rationally connected. 
So we can find a very free curve $l_1$ through $P_1$ on this fiber, and consider the curve
$C\cup l_1$. We repeat the process for $k\gg 1$ general points in $C$,  
obtaining a reducible curve $C\cup l_1\cup\dots \cup l_k$.
Then we deform this reducible curve
into a single irreducible curve $C'$. 
Roughly speaking, we can make the space of deformations of $C'$
as big as we want by taking $k$ large enough.


We introduce some concepts and gather some results that will 
be used in the proof of Theorem~\ref{ghs}.

\subsection{Stable Maps and Hurwitz schemes}\label{moduli_spaces}

\begin{defn}[Stable maps]
Let $X$ be a complex projective variety.
We say that a morphism $h:C\to X$ is a stable map if 
\begin{enumerate}
  \item $C$ is a projective connected curve with at worst nodes as singularities,
  \item if  $C'$ is an irreducible component of $C$ that has arithmetic genus $0$ and is 
    contracted by $h$, then  $C'$ contains at least three nodes of $C$, and
  \item if  $C'$ is an irreducible component of $C$ that has arithmetic genus $1$ and is
    contracted by $h$, then  $C'$ contains at least one node of $C$.
\end{enumerate}
\end{defn}

Fix a nonnegative integer $g$ and a class $\alpha$ in $N_1(X)$.
There is a projective scheme $\overline M_{g}(X,\alpha)$
that is a coarse moduli space for all stable maps $h:C\to X$
such that $C$ has arithmetic genus $g$, and $h_*[C]=\alpha$
(see \cite{fulton_pandharipande}).

Now let $f:X\to B$ be a morphism onto a smooth curve, and 
set $d=f_*\alpha$. There exists a natural morphism
$$
\bar f: \overline M_{g}(X,\alpha)\to \overline M_{g}(B,d)
$$ 
obtained by composing a stable map $h:C\to X$
with $f$ and collapsing the irreducible components of $C$ that make $f\circ h$ unstable
(see \cite[Theorem 3.6]{behrend_fantechi}).

\begin{defn}[Space of deformations]
Let $C\subset X$ be a projective connected curve with at worst nodes as singularities.
Let $g$ be the arithmetic genus of $C$, and $\alpha$ the class of $C$ in $N_1(X)$. 
Let $h:C\to X$ be the inclusion morphism. 
Then $h$ is a stable map, and we denote by $[C]$
the point of $\overline M_{g}(X,\alpha)$ parametrizing $h$.

We define the space of deformations of $C$ in $X$, $\defm(C,X)$, to be union of
the irreducible components of $\overline M_{g}(X,\alpha)$ containing $[C]$. 
\end{defn}

Let $C\subset X$ be a projective connected curve with at worst nodes as singularities.
There is a nice description of the Zariski tangent space of $\defm(C,X)$ at the point $[C]$:
$$
T_{[C]}\defm(C,X)\cong H^0(C,N_C),
$$ 
where $N_C$ denotes the normal bundle of $C$ in $X$.
If $H^1(C,N_C)=0$, then $\defm(C,X)$ is smooth at $[C]$.
This means that every first order deformation $\alpha \in H^0(C,N_C)$
comes from a $1$-parameter family of deformations of $C$ in $X$.

\begin{defn}[Hurwitz Schemes]
Fix a positive integer $d$, a nonnegative integer $g$, and set $\delta=2g+2d-2$.

Let $C$ be a smooth complex projective irreducible curve of genus $g$. Let 
$h:C\to \p^1$ be a $d$-sheeted branched covering, and $B\subset \p^1$ the branch divisor of $h$.
By the Riemann--Hurwitz formula, $\deg B = 2g+2d-2 = \delta$.
We say that $h$ is simply branched if $B$ is reduced. 

We define the Hurwitz scheme $H^{d,\delta}$ to be the open subscheme of $\overline M_{g}(\p^1,d)$
parametrizing $d$-sheeted simply branched coverings of $\p^1$. 
We denote by $\overline H^{d,\delta}$ the closure of $H^{d,\delta}$ in $\overline M_{g}(\p^1,d)$.
\end{defn}

It turns out that the Hurwitz scheme $H^{d,\delta}$ is smooth and irreducible 
(see \cite[Proposition 1.5]{fulton_on_hurwitz}). Its closure $\overline H^{d,\delta}$
contains all $d$-sheeted branched coverings from smooth curves of genus $g$. 
It also contains points parametrizing morphisms $h:D\to \p^1$ with a section. Indeed, let $h':D'\to \p^1$ 
consist of $d$ disjoint copies of $\p^1$, each one mapping isomorphically onto $\p^1$. 
Then $h:D\to \p^1$ can be obtained from $h'$ by identifying suitable $d+g-1$ pairs of points in $D'$, each pair
lying in a different fiber of $h'$.

Consider the branch divisor map 
$$
\br: \ H^{d,\delta}\ \to \ \sym^{\delta}(\p^1),
$$
which assigns to each branched covering its branch divisor.
The map $\br$ is an \'etale cover of its image.
By \cite{fantechi_pandharipande}, $\br$ extends to a morphism
$$
\overline \br: \ \overline H^{d,\delta}\ \to \ \sym^{\delta}(\p^1).
$$

\subsection{The proof of Theorem~\ref{ghs}}\label{proof_GHS}

Now we explain the proof of Theorem~\ref{ghs} in the case when $B\cong \p^1$. 
At the end of the section we explain how the proof can be extended to arbitrary curves $B$.
We assume that $\dim X\geq 4$. 
Otherwise replace $X$ with $X'=X\times \p^N$, and $f$ with $f'$, the composition of $f$ with 
the projection onto $X$. 
Any section of $f'$ yields a section of $f$.

Let $C\subset X$ be a smooth projective irreducible curve for which  $f_C=f|_C$ is a finite morphism
(we say that $C$ is a multisection of $f$).
Let $g$ be the genus of $C$, let $d$ be the degree of $f_C$, and set $\delta=2g+2d-2$.
Let $R\subset C$ be the ramification divisor of $f_C$, and
$B\subset \p^1$ the branch divisor of $f_C$.

By composing the morphisms $\bar f$ and $\overline \br$ defined above, we obtain the following branch morphism:
$$
\Br: \  \defm(C,X)\ \overset{\bar f}{\longrightarrow} \ \overline H^{d,\delta}\ 
\overset{\overline \br}{\longrightarrow} \ \sym^{\delta}(\p^1).
$$

If this morphism is surjective, then we are done. 
Indeed, if $\Br$ is surjective, then so is $\bar f$ (as $\overline \br$ is finite).
As we saw above, $\overline H^{d,\delta}$
contains a point $[h]$ parametrizing a morphism $h:D\to \p^1$ with a section.
Any point of $\defm(C,X)$ sent to $[h]$ by $\bar f$
yields a curve on $X$ 
containing a section of $f$ as an irreducible component.

Assume that $f_C$ is simply branched (so that $\sym^{\delta}(\p^1)$ is smooth at $[B]$)
and $\defm(C,X)$ is smooth at $[C]$.
In order to show that $\Br$ is surjective, it suffices to show that the derivative map 
$$
d\Br_{[C]}:T_{[C]} \defm(C,X)\to T_{[B]}\sym^{\delta}(\p^1)
$$ 
is surjective.

Let us describe this map explicitly.
First we note that there is an isomorphism $T_{[C]}\defm(C,X)\cong H^0(C,N_C)$.
Let $\ev:H^0(C,N_C)\to (N_C)|_R$ denote the map obtained by evaluating sections 
at the ramification points of $f_C$.
Since the derivative map $df_C$ is zero at each ramification point of $f_C$,
the derivative of $f$ induces a map $df:(N_C)|_R\to (T_{\p^1})|_B\cong T_{[B]}\sym^{\delta}(\p^1)$.
We then have $d\Br_{[C]}=df\circ \ev$:
\[
\xymatrix{
H^0(C,N_C) \ar[rr]^{d\Br_{[C]}} \ar[rd]_{\ev} & & (T_{\p^1})|_B. \\
& (N_C)|_R \ar[ru]_{df} 
}
\]

From this discussion we conclude that to prove Theorem~\ref{ghs}
it suffices to find a multisection $C\subset X$ satisfying the following conditions:

\begin{condition}\label{conditions}
\  
  \begin{enumerate}
    \item $f_C$ is simply branched, 
    \item $\defm(C,X)$ is smooth at $[C]$,
    \item the map $\ev:H^0(C,N_C)\to (N_C)|_R$ is surjective, and
    \item the map $df:(N_C)|_R\to (T_{\p^1})|_B$ is surjective.
  \end{enumerate}
\end{condition}

In order to achive this, we will make use of two constructions.
The first construction allows us to, starting with any multisection $C$,  
deform it into a new multisection $C'$ for which the evaluation map $\ev'$ is surjective.
As a bonus, we also get that  $\defm(C',X)$ is smooth at $[C']$.
The second construction allows us to, starting with any multisection $C$, deform it and 
degenerate it so as to produce a new multisection $C'$ for which 
$df$  is surjective.

\begin{construction}\label{construction_1}
Let $C\subset X$ be a smooth multisection of $f$.
Let $R\subset C$ be the ramification divisor of $f_C$. 

Let $P_1,\dots ,P_k\in C$ be general points.  
By assumption, the fibers of $f$ containing each $P_i$ are rationally connected. 
So we pick $l_1,\dots,l_k\subset X$ general smooth very free curves on fibers of $f$, each $l_i$ meeting $C$
at $P_i$ only. 
Such curves exist by Remark~\ref{remark_v_free} and Proposition~\ref{avoiding_codim_2}.
(Here we need the assumption that the fibers of $f$ have dimension at least $3$.)
We call each such $l_i$ a \emph{very free tail on $C$}. 
We show below that, if $k$ is large enough, then we can deform 
$C\cup l_1\cup\dots \cup l_k$ into a single smooth irreducible curve $C'$.
Moreover, by increasing $k$ if necessary, 
we have $H^1(C',N_{C'}(-R'))=0$, where  
$R'\subset C'$ is the ramification divisor of $f_{C'}:C'\to \p^1$.
(Roughly speaking, the more very free tails we attach to $C$, the more positive the normal bundle 
of the resulting curve becomes.)

Let $\ev':H^0(C',N_{C'})\to (N_{C'})|_{R'}$ be the evaluation map. 
From the exact sequence 
$$
0 \ \to \  N_{C'}(-R')\ \to \ N_{C'} \  \to \  (N_{C'})|_{R'}  \ \to \ 0, 
$$
we see that $\coker(\ev')\subset H^1(C',N_{C'}(-R'))=0$, and thus
$\ev'$ is surjective.

Moreover, $H^1(C',N_{C'}(-R'))=0$ implies that $H^1(C',N_{C'})=0$.
Hence the space of deformations  $\defm(C',X)$ is smooth at $[C']$.



\begin{figure}[!htb]
\centering
\psfrag{X}{$X$}
\psfrag{C}{$C$}
\psfrag{P}{$\p^1$}
\psfrag{f}{$f$}
\psfrag{l1}{$l_1$}
\psfrag{l2}{$l_2$}
\psfrag{l3}{$l_3$}
\psfrag{lk}{$l_k$}
\begin{minipage}[b]{0.45\linewidth}
\includegraphics[width=\linewidth]{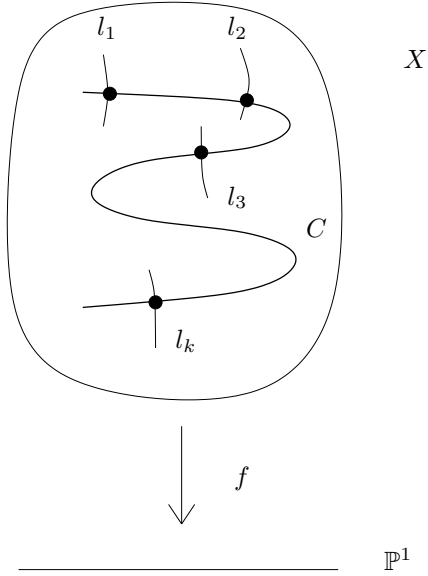}
\caption{Very free tails} 
\label{tails}
\end{minipage} \hfill
\end{figure} 


\end{construction}

Next we prove that Construction~\ref{construction_1} in fact works to produce a smooth multisection $C'$
satisfying $H^1(C',N_{C'}(-R'))=0$, and hence
Conditions~\ref{conditions} (2) and (3). 
The proof we give here is slightly different from the one in \cite{GHS}, and follows an 
argument in \cite{kollar_specialization}.
The reader may want to skip it in a first reading.

\begin{proof}[Proof that Construction~\ref{construction_1} works]

Let $C$ be a smooth multisection of $f$.
Let $g$ be the genus of $C$, let $d$ be the degree of $f_C$, and set $\delta=2g+2d-2$.
Let $D$ be any divisor on $C$ of degree $\delta + g$.
Our aim is to attach to $C$ very free tails $l_1,\dots,l_k$ disjoint from $D$
so that
\begin{equation}\label{eq:vanishing}
H^1(C\cup l_1\cup \dots \cup l_k, N_{C\cup l_1\cup \dots \cup l_k}(-D))=0.
\end{equation}

\medskip

In order to produce such a curve $C\cup l_1\cup \dots \cup l_k$,
we use induction on the number of very free tails on $C$.

Let $l_1,\dots, l_m$ be very free tails on $C$ disjoint from $D$.
Consider the curve 
$C_m=C\cup l_1\cup \dots \cup l_m$.
We have exact sequences
$$
0\ \to \ N_{l_i}(-P_i)\ \to \ (N_{C_m})|_{l_i}(-P_i)\ \to \ Q_i \ \to \ 0,
$$
where $Q_i$ is a torsion sheaf supported at $P_i=C\cap l_i$, and 
$$
0\ \to \ \bigoplus_{i=1}^m (N_{C_m})|_{l_i}(-P_i)\ \to \ N_{C_m}(-D)\ \to \ (N_{C_m})|_C(-D) \ \to \ 0.
$$

From these exact sequences and the assumption that each $l_i$ is a very curve on a fiber of $f$,
it follows that $H^1(C_m,N_{C_m}(-D)) \cong H^1(C,(N_{C_m})|_C (-D))$.
By Serre duality, the latter is dual to 
$\Hom (\omega_C^{-1}(-D), (N_{C_m}^{\vee})|_C)$.
Suppose that this does not vanish. Pick a nonzero element
$\alpha \in \Hom (\omega_C^{-1}(-D), (N_{C_m}^{\vee})|_C)$.

Let $P_{m+1}\in C$ be a general point. 
Then $P_{m+1}\notin l_1\cup \dots \cup l_m$, and 
the fiber of $f$ through $P_{m+1}$ is rationally connected.
Moreover, the composition
$$
\bar \alpha:\ \omega_C^{-1}(-D)\ \overset{\alpha}{\longrightarrow} (N_{C_m}^{\vee})|_C \
\longrightarrow \ (N_C^{\vee})|_{P_{m+1}}
$$
has rank $1$.
Let $\phi:(N_C)|_{P_{m+1}}\to \c$ span the image of $\bar \alpha$, and let 
$\xi\in (N_C)|_{P_{m+1}}$ be a general normal vector at $P_{m+1}$, 
so that $\xi\notin \ker \phi$.
By Remark~\ref{remark_v_free}(6), 
there exists a smooth very free curve $l_{m+1}$ on the fiber of $f$ through $P_{m+1}$
with normal direction $\xi$ at $P_{m+1}$.
We also require that $l_{m+1}$ meets $C$ at $P_{m+1}$ only. 
Consider the curve $C_{m+1}=C\cup l_1\cup \dots \cup l_{m+1}$.
We have the exact sequence
$$
0\ \to \ (N_{C_{m+1}}^{\vee})|_C\ \to \ (N_{C_m}^{\vee})|_C\ \overset{\beta}{\rightarrow} \c \xi ^{\vee} \ \to \ 0.
$$

Let us look again at $\alpha \in \Hom (\omega_C^{-1}(-D), (N_{C_m}^{\vee})|_C)$
and $\phi\in (N_{C}^{\vee})|_{P_{m+1}}$ spanning the image of 
$\bar \alpha:\omega_C^{-1}(-D)\to (N_C^{\vee})|_{P_{m+1}}$.
By construction, $\xi$ does not lie in the kernel of $\phi$.
Hence $\beta(\alpha)\neq 0$, and $\alpha$
does not come from an element in  $\Hom (\omega_C^{-1}(-D), (N_{C_{m+1}}^{\vee})|_C)$.
Thus 
\begin{equation*}
\begin{split}
h^1(C_m,N_{C_m}(-D)) & = \dim \Hom (\omega_C^{-1}(-D), (N_{C_{m}}^{\vee})|_C) \\
& < \dim \Hom (\omega_C^{-1}(-D), (N_{C_{m+1}}^{\vee})|_C) \\
& = h^1(C_{m+1},N_{C_{m+1}}(-D)).
\end{split}
\end{equation*}

So, after finitely many steps, we obtain a curve $C_k=C\cup l_1\cup \dots \cup l_k$ 
satisfying \eqref{eq:vanishing} above. 

\medskip

Now we are pretty much done.
It is easy to see that \eqref{eq:vanishing} implies that $N_{C_k}$ is generated by
global sections.
So we can find a global section $\alpha \in H^0(C_k,N_{C_k})$ such that, 
for every  $i \in \{1,\dots,k\}$, the restriction 
$\alpha|_{P_i}\in (N_{C_k})|_{P_i}$ does not come from $(N_C)|_{P_i}$.
Since $H^1(C_k,N_{C_k})=0$, $\alpha$ comes from a global deformation of $C_k$
in $X$ smoothing its nodes at each $P_i$.
Let $C'$ be a general such deformation, and $D'\subset C'$ the corresponding deformation of $D$.
By Semicontinuity, $H^1(C',N_{C'}(-D'))=0$, which implies, by Riemann--Roch, that 
$H^1(C',N_{C'}(-R'))=0$, where $R'$ is the ramification divisor of 
$f_{C'}:C'\to \p^1$.
\end{proof}

Before we explain the second construction, let us investigate conditions under which 
we have a smooth multisection $C$ for which the corresponding map 
$df:(N_C)|_R\to (T_{\p^1})|_B$ is surjective.

First notice that $df$ is surjective provided that $C$ is contained in the smooth locus of $f$.
Denote by $\sing(f)$ the singular locus of $f$.
Suppose now that $C$ only meets irreducible components of $\sing(f)$ that 
have codimension at least $2$ in $X$.
By applying Construction~\ref{construction_1} if necessary, we may assume that 
$H^1(C,N_{C}(-P))=0$ for every $P\in C$.
In this case one can show that the general deformation of $C$ that does not meet $\sing(f)$ at all.
So we only have a problem when $C$ meets some irreducible component of $\sing(f)$ 
of codimension $1$ in $X$. 
Such irreducible component is necessarily a multiple component of a fiber of $f$.
We take care of this case in Construction~\ref{construction_2} below.

\begin{construction}\label{construction_2}
Let $C\subset X$ be a smooth multisection of $f$.
Our aim is to produce a new multisection $C'$ that does not meet 
any multiple fiber of $f$.
Let $g$ be the genus of $C$, let $d$ be the degree of $f_C$, and set $\delta=2g+2d-2$.
Let $R\subset C$ be the ramification divisor of $f_C$, 
and $B\subset \p^1$ the branch divisor of $f_C$.
As we noted above, we may assume that $C$ only meets irreducible components 
of $\sing(f)$ of codimension $1$ in $X$.

We introduce the monodromy homomorphism
associated to the branched covering $f_C:C\to \p^1$.
First choose a base point $b\in \p^1\setminus B$, and fix an identification
$\Phi:f_C^{-1}\{b\}\overset{\cong}{\longrightarrow} \{1,\dots, d\}$.
Then define the monodromy homomorphism
$$
\tilde \Phi: \pi_1(\p^1\setminus B,b)\to S_d
$$
by associating to a loop $\gamma\in \pi_1(\p^1\setminus B,b)$ 
the permutation of $\{1,\dots, d\}$ obtained by analytic continuation 
along $\gamma$.

Let $b_1,\dots, b_m \in \p^1$ be the branch points of $f_C$ whose fibers 
contain multiple components. 
Suppose we can find disjoint discs $\Delta_i\subset \p^1$, $1\leq i\leq m$, such that
$b_i\in \Delta_i$, $b\in \partial \Delta_i\subset \p^1\setminus B$, and
\begin{equation}\label{eq:trivial_monodromy}
\text{{\it the monodromy around the boundary }}\partial \Delta_i\text{ {\it is trivial.}}
\tag{*}
\end{equation}
Notice that each $\Delta_i$ must contain other branch points of $f_C$
distinct from $b_1,\dots, b_m$.
Let $s_1,\dots, s_l$ be the branch points of $f_C$ that do not lie in any 
$\Delta_i$.



\begin{figure}[!htb]
\centering
\psfrag{b}{$b$}
\psfrag{b1}{$b_1$}
\psfrag{P}{$\p^1$}
\psfrag{b2}{$b_2$}
\psfrag{bm}{$b_m$}
\psfrag{g1}{{\bf $\partial \Delta_1$}}
\psfrag{g2}{{\bf $\partial \Delta_2$}}
\psfrag{gm}{{\bf $\partial \Delta_m$}}
\begin{minipage}[b]{0.5\linewidth}
\includegraphics[width=\linewidth]{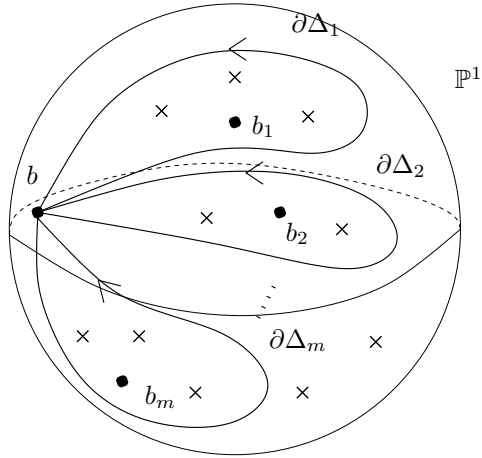}
\caption{Loops with trivial monodromy} 
\label{monodromy}
\end{minipage} \hfill
\end{figure} 


The idea now is to keep the ``bad'' branch points $b_i$ fixed and move
the other branch points in $\Delta_i$ toward $b_i$, also keeping the 
branch points $s_j$ fixed. Then we lift this deformation of $B$ in
$\sym^{\delta}(\p^1)$ to a deformation of $C$ in $X$ 
and look at what happens in the limit.

Let $Z$ be the following analytic subvariety of $\sym^{\delta}(\p^1)$.
$$
Z \ = \ \sum_{i=1}^m m_{b_i}b_i + \sum_{i=1}^m \tilde m_i \Delta_i  + \sum_{j=1}^l m_{s_j}s_j ,  
$$
where $m_{b_i}$ is the multiplicity of $b_i$ in $B$,
$m_{s_j}$ is the multiplicity of $s_j$ in $B$ and
$\tilde m_i$ is the number of branch points of $f_C$ in $\Delta_i\setminus \{b_i\}$,
counted with multiplicity.

Let $0,t_0\in \Delta$ be two distinct points in the unit disc.
Let $B_t$, $t\in \Delta$, be an analytic arc in $Z$ such that 
\begin{enumerate}
  \item $B_{t_0}=B$, and
  \item $B_0=\sum_{i=1}^m(m_{b_i}+\tilde m_i)b_i+\sum_{j=1}^l m_{s_j}s_j$.
\end{enumerate}



\begin{figure}[!htb]
\centering
\psfrag{b}{$b$}
\psfrag{b1}{$b_1$}
\psfrag{P}{$\p^1$}
\psfrag{b2}{$b_2$}
\psfrag{bm}{$b_m$}
\psfrag{g1}{$\gamma_1$}
\psfrag{g2}{$\gamma_2$}
\psfrag{gm}{$\gamma_m$}
\begin{minipage}[b]{0.5\linewidth}
\includegraphics[width=\linewidth]{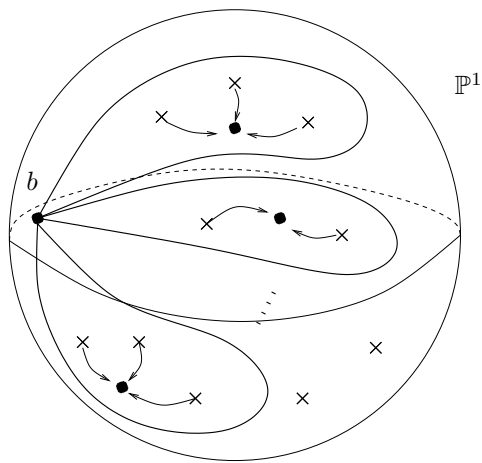}
\caption{The analytic arc $B_t\subset \sym^{\delta}(\p^1)$} 
\label{monodromy2}
\end{minipage} \hfill
\end{figure} 


By applying Construction~\ref{construction_1} to $C$ if necessary, we may assume that 
the map $\ev:H^0(C,N_C)\to (N_C)|_R$ is surjective.
Moreover, the map $df:(N_C)|_P\to (T_{\p^1})|_{f(P)}$ is surjective for every 
$P\in R$ such that $f(P)\neq b_i$ for $1\leq i\leq m$. 
Therefore $\Br: \defm(C,X)\to \sym^{\delta}(\p^1)$ dominates $Z$, and
the restriction $\Br|_{\Br^{-1}Z}:\Br^{-1}Z\to Z$ is smooth at $[C]\in \defm(C,X)$.
So we can find an analytic arc $\gamma_t\subset \defm(C,X)$ lifting $B_t$ 
with $\gamma_{t_0}=[C]$.

Let us look at the limit $\gamma_0\in \defm(C,X)$. 
It corresponds to a stable map $h_0:C_0\to X$.
Let $C_0'$ be an irreducible component of $C_0$ that dominates
$\p^1$ via $f\circ h_0$, and let $f_{C_0'}:C_0'\to \p^1$ 
be the induced branched covering.
By construction,
\begin{enumerate}
  \item the only possible branch point of $f_{C_0'}$ in $\Delta_i$ is $b_i$
    (indeed $\Br(\gamma_0)=B_0$), and
  \item the monodromy around the boundary $\partial \Delta_i$ is trivial.
\end{enumerate} 
These two conditions together imply that $f_{C_0'}$ is unramified over $\Delta_i$
for $1\leq i\leq m$, and it can only ramify over the $s_j$
with multiplicity at most $m_{s_j}$. 
Hence $C'=h_0(C_0')$ does not meet any multiple fiber of $f$.
Notice that if $f_C$ is simply branched over each $s_j$, then
$f_{C_0'}$ can be at most simply branched over each $s_j$.

After adding some very free tails and deforming if necessary,
we may assume that $C'$ is smooth and is contained in the smooth locus of $f$.
So it satisfies Condition~\ref{conditions} (4).
Notice that if $f_C$ is simply branched over each $s_j$, then
$C'$ also satisfies Condition~\ref{conditions} (1).

Now we have to find disjoint discs $\Delta_i$, $1\leq i\leq m$, with $b_i\in \Delta_i$ and 
satisfying condition \eqref{eq:trivial_monodromy} above.
In order to do so, we need that $f_C$ has enough ``good'' branch points
with suitable monodromy to cancel out the monodromy around each $b_i$.
If this is not the case, then we will need to create new branch points with
preassigned monodromy. Next we explain how this can be done.
Notice that, since $S_d$ is generated by transpositions, it is enough to assign
transpositions as monodromies.

Let $t\in \p^1\setminus (B\cup \{b\})$ be a general point.
Let $i,j$ be distinct elements of $\{1,\dots,d\}\cong f_C^{-1}\{b\}$.
Fix a (real) path in $\p^1\setminus B$ joining $b$ and $t$.
By analytic continuation along this path, $i$ and $j$ determine
two points $x$ and $y$ in $f_C^{-1}\{t\}$.
The fiber of $f$ over $t$ is rationally connected.
So, by Remmark~\ref{remark_v_free}(5), there exists a smooth very free curve $l$ 
on this fiber connecting $x$ and $y$.
We may also require that $l$ meets $C$ at $x$ and $y$ only.
After attaching some very free tails to $C$ if necessary, we are able to 
deform $C\cup l$ into a smooth irreducible curve $\tilde C$.
(We skip the proof as it is very similar to the one given in 
Construction~\ref{construction_1}.)
The branched covering $f_{\tilde C}:\tilde C\to \p^1$ has the same degree as $f_C$.
Its branch divisor is a deformation of $B$ plus two new simple branch points.
The monodromy (of suitable loops) around each one of the new branch points is precisely
$(ij)$.

After creating enough new branch points and replacing $C$ with the resulting curve,
we can find disjoint discs $\Delta_i$ with $b_i\in \Delta_i$ and satisfying condition~\eqref{eq:trivial_monodromy}.
Then we degenerate $C$ as explained above, obtaining a smooth multisection $C'$
of $f$ satisfying Condition~\ref{conditions} (4). 
\end{construction}

Next we explain how to obtain a smooth multisection of $f$ satisfying Conditions~\ref{conditions}(1)--(4).

First embed $X$ in some projective space.
By intersecting $X$ with $\dim X-1$ general hyperplane sections, we get a smooth multisection $\tilde C$ of $f$.
Dimension counts show that $f_{\tilde C}$ is at most simply branched over the points of $\p^1$
whose fibers do not contain multiple components.

Now apply Construction~\ref{construction_2} to $\tilde C$: 
create branch points with assigned monodromy, make the normal bundle of the resulting curve very positive
by attaching very free tails and deforming, and degenerate it by colliding some of the branch points 
and keeping the other ones fixed. 
This yields a smooth multisection $C'$ of $f$ satisfying Conditions~\ref{conditions} (1) and (4).

Finally, apply Construction~\ref{construction_1} to $C'$:
attach many very free tails to $C'$ and deform the resulting reducible curve into a smooth multisection $C$.
Notice that $C$ still satisfies Conditions~\ref{conditions} (1) and (4), as these are open conditions.
It also satisfies Conditions~\ref{conditions} (2) and (3).
This ends the proof of Theorem~\ref{ghs}.

\subsection{Exending the proof to arbitrary curves {\it B}}

Among the ingredients used in the proof of  Theorem~\ref{ghs} above, 
the only one that relies on the assumption that $B\cong \p^1$ is the irreducibility of the 
Hurwitz scheme $H^{d,\delta}$.

Let $B$ be a smooth complex projective irreducible curve of genus $h$.
Fix a positive integer $d$ and a nonnegative integer $g$ such that $\delta=2g-2d(h-1)-2\geq 0$.
Let $G$ be a subgroup of the symmetric group $S_d$.
We define the Hurwitz scheme $H^{d,\delta}_G(B)$ to be the open subscheme of $\overline M_{g}(B,d)$
parametrizing $d$-sheeted simply branched coverings of $B$ 
whose monodromy group is conjugate to $G$ in $S_d$.

The Hurwitz scheme $H^{d,\delta}_G(B)$ is not irreducible in general.
It is irreducible, however, if the following two conditions are satisfied.
\begin{enumerate}
  \item $\delta \geq 2d$, and
  \item $G=S_d$.
\end{enumerate}
(See \cite{GHS_Hurwitz} for a proof.)

In order to extend the proof of Theorem~\ref{ghs} to arbitrary curves $B$,
this is what we do.
Start with a smooth multisection $C$ obtained by 
intersecting $\dim X-1$ general very ample divisors on $X$. 
Consider the branched covering $f_C:C\to B$.
If $f_C$ contains a large number of branch points (more than $2\deg f_C$), and 
if its monodromy group is the full symmetric group $S_{\deg f_C}$, then 
the argument given above goes through.
Otherwise, create new branch points with assigned monodromy, as 
explained in Construction~\ref{construction_2}, so that 
the resulting curve $\tilde C$ satisfies these two conditions.
Then replace $C$ with $\tilde C$ and proceed with the argument.


\section{Further results and open problems}\label{section:final}

\subsection{Rationally connected varieties in positive characteristic}

Fix an algebraically closed field $\k$.
Let $X$ be a smooth projective variety over $\k$.
If $\k$ has positive characteristic, 
the existence of rational curves through two general points does not imply 
the existence of a very free curve on $X$. 
So we have two different notions.

We say that $X$ is rationally connected if  the following holds. 
There exists an irreducible component $H$ of $\hilb(X)$ whose general member parametrizes 
a rational curve and for which the double-evaluation map
$\cU\times_H \cU \overset{\ev^{(2)}}{\longrightarrow}  X \times X$
is dominant. 
Here $\cU$ is the universal family over $H$, and $\ev:\cU\to X$ 
the evaluation map.
If $\k$ is uncountable, this condition is equivalent to 
the existence of rational curves through two general points.

We say that $X$ is \emph{separably rationally connected} if there 
exists an irreducible component of $\hilb(X)$ as above 
for which the double-evaluation map $\ev^{(2)}$ is \emph{separable} and dominant.
This condition is equivalent to the existence of a very free curve on $X$.

The following generalization of Theorem~\ref{ghs} was proved by de Jong and Starr in \cite{dejong_starr}.

\begin{thm}[\cite{dejong_starr}] \label{ghs_over_k} 
Let $X$ be a smooth projective variety over $\k$.
Let $f:X\to B$ be a surjective morphism onto a smooth curve, and assume that the 
general fiber of $f$ is separably rationally connected. 
Then $f$ has a section. 
\end{thm}

\subsection{A converse of Theorem~\ref{ghs}}

Let $X$ and $Y$ be complex varieties and $f:X\to Y$ a proper morphism.
Suppose that the restriction $f|_{f^{-1}C}:f^{-1}C\to C$
has a section for every curve $C\subset Y$.
The strict converse of  Theorem~\ref{ghs} would assert that the general fiber of $f$ is 
rationally connected. 
This is clearly false. 
Take, for instance, $X=Y\times Y'$, where $Y'$ is not rationally connected, and $f$ the 
projection onto the first factor.
However, something weaker holds:
there exists a subvariety $Z\subset X$ for which the restriction $f|_Z:Z\to Y$ is 
dominant with rationally connected general fiber.
In \cite{GHS_converse} such $Z$ is called a \emph{pseudosection} of $f$.
Recall that a point is rationally connected,
so a section of $f$ is a pseudosection.

In fact, the result proved in \cite{GHS_converse} is stronger than this.
It says that, instead of checking whether $f|_{f^{-1}C}:f^{-1}C\to C$
has a section for every curve $C\subset Y$, it suffices to check it for a
fixed bounded family of test curves.

\begin{thm}[\cite{GHS_converse}]
Let $Y$ be a complex variety. Fix a positive integer $d$.
Then there exists a bounded family $H_d$ of irreducible curves on $Y$ satisfying the following property.
If $f:X\to Y$ is a proper morphism of relative dimension at most $d$, 
and if $f|_{f^{-1}C}:f^{-1}C\to C$ has a section for a very general curve $C$ parametrized by $H_d$,
then there exists a subvariety $Z\subset X$ for which the restriction $f|_Z:Z\to Y$ is 
dominant with rationally connected general fiber.
\end{thm}

\subsection{Producing sections with preassigned data}

Let $X$ be a smooth complex projective variety and $f:X\to B$ a surjective morphism onto a smooth curve. 
Suppose that the general fiber of $f$ is rationally connected.
Once Theorem~\ref{ghs} gives us a section of $f$ to start with, 
it is easy to produce new sections passing through finitely many preassigned points
in smooth fibers of $f$.
The basic idea is the following. 
Let $C\subset X$ be a section of $f$, $X_{b_1},\dots,X_{b_k}$ distinct smooth fibers of $f$,
and $x_i\in X_{b_i}$, $1\leq i \leq k$.
Since each $X_{b_i}$ is rationally connected,
there are very free rational curves $l_i$ on $X_{b_i}$ connecting $x_i$ and $C\cap X_{b_i}$.
By attaching very free tails to $C\cup l_1 \cup \dots \cup l_k$, 
we are able to deform the resulting curve keeping the $x_i$ fixed. 
In this way we produce the required section 
(see \cite[IV.6.10]{kollar} for details).

In \cite{hassett_tschinkel}, Hassett and Tschinkel improved this argument and showed that we 
can in fact prescribe finitely many terms of the Taylor expansion of the section 
$\sigma: B\to X$ at the points $x_i$. More precisely:

\begin{thm}[\cite{hassett_tschinkel}]\label{hassett_tschinkel}
Let $X$ be a smooth complex projective variety and $f:X\to B$ a surjective morphism onto a smooth curve. 
Suppose that the general fiber of $f$ is rationally connected.
Let $b_1, \dots, b_k\in B$ be distinct points whose fibers are smooth.
Let $N$ be a positive integer and for each $i\in \{1,\dots,k\}$ 
let $\sigma_i$ be a section of the induced morphism
$$
X\times_B \spec(\o_{B,b_i}/\m_{B,b_i}^N) \ \rightarrow \ \spec(\o_{B,b_i}/\m_{B,b_i}^N).
$$
Then $f$ has a section agreeing with all the $\sigma_i$.
\end{thm}

They also conjectured that one may relax the assumption that the fibers $X_{b_i}$ are 
smooth, and instead require that each $\sigma_i(b_i)$ is a smooth point of $X_{b_i}$.

\subsection{Open problems}

Here we only discuss a couple of open problems concerning rationally connected varieties
in the context presented in section~\ref{intro}. 
We refer to \cite{kollar_simplest_alg_var} for a more complete list of open problems in the area.

\medskip

As we noted in section~\ref{intro}, there is a conjectural numerical characterization of
complex rationally connected verieties.

\begin{conj}\label{conj_RCV}
Let $X$ be a smooth complex projective variety.
Then $X$ is rationally connected if and only if $H^0(X,(\Omega_X^1)^{\otimes m})=0$ for every $m\geq 1$.
\end{conj}

The ``only if'' part follows very easily from the definition, as we showed in section~\ref{section:RCV}.
The ``if'' part is a very hard problem. 
Theorem~\ref{ghs} reduces it to a similar conjecture about uniruled varieties (see \cite{GHS} for details).

\begin{conj}\label{conj_uniruled}
Let $X$ be a smooth complex projective variety.
Then $X$ is uniruled if and only if $H^0(X,\omega_X^{\otimes m})=0$ for every $m\geq 1$.
\end{conj}

As in Conjecture~\ref{conj_RCV}, the ``only if'' part of Conjecture~\ref{conj_uniruled} follows 
easily from the definition.

\medskip

There is another class of varieties lying between the class of rational varieties and the class of 
rationally connected varieties.

\begin{defn}
Let  $X$ be a smooth projective variety of dimension $n$.
We say that $X$ is \emph{unirational} if there is a dominant map $\p^n\map X$.
\end{defn}

The class of unirational varieties is strictly bigger than the class of rational varieties.
Smooth cubic $3$-folds, for instance, are unirational but not rational. 
Unirational varieties are clearly rationally connected.
On the other hand, it is not known whether the classes of unirational varieties and 
rationally connected varieties are in fact distinct.

\begin{problem}
Find examples of rationally connected varieties that are not unirational. 
\end{problem}

\bigskip

\noindent {\it Acknowledgments.\ } These notes are based on the lecture that I gave at Snowbird, 2004.
Some ideas for the presentation of the proof of Theorem~\ref{ghs} were inspired by a lecture 
given by O. Debarre at the University of Warwick in 2002.
I would like to thank J. Ellenberg, S. Grushevsky, S. Kov\'acs and J. Starr for many comments and suggestions.

\bibliographystyle{amsalpha}
\bibliography{GHS_bib}

\providecommand{\bysame}{\leavevmode\hbox to3em{\hrulefill}\thinspace}
\providecommand{\MR}{\relax\ifhmode\unskip\space\fi MR }
\providecommand{\MRhref}[2]{%
  \href{http://www.ams.org/mathscinet-getitem?mr=#1}{#2}
}
\providecommand{\href}[2]{#2}
\begin{thebibliography}{KMM92b}

\bibitem[BF97]{behrend_fantechi}
K.~Behrend and B.~Fantechi, \emph{The intrinsic normal cone}, Invent. Math.
  \textbf{128} (1997), no.~1, 45--88.

\bibitem[Cam92]{campana}
F.~Campana, \emph{Connexit\'e rationnelle des vari\'et\'es de {F}ano}, Ann.
  Sci. \'Ecole Norm. Sup. (4) \textbf{25} (1992), no.~5, 539--545.

\bibitem[CG72]{clemens_griffiths}
H.~Clemens and P.~Griffiths, \emph{The intermediate {J}acobian of the cubic
  threefold}, Ann. of Math. (2) \textbf{95} (1972), 281--356.

\bibitem[Deb01]{debarre}
O.~Debarre, \emph{Higher-dimensional algebraic geometry}, Universitext,
  Springer-Verlag, New York, 2001.

\bibitem[dJS03]{dejong_starr}
A.~J. de~Jong and J.~Starr, \emph{Every rationally connected variety over the
  function field of a curve has a rational point}, Amer. J. Math. \textbf{125}
  (2003), no.~3, 567--580.

\bibitem[FP97]{fulton_pandharipande}
W.~Fulton and R.~Pandharipande, \emph{Notes on stable maps and quantum
  cohomology}, Algebraic geometry---Santa Cruz 1995, Proc. Sympos. Pure Math.,
  vol.~62, Amer. Math. Soc., Providence, RI, 1997, pp.~45--96.

\bibitem[FP02]{fantechi_pandharipande}
B.~Fantechi and R.~Pandharipande, \emph{Stable maps and branch divisors},
  Compositio Math. \textbf{130} (2002), no.~3, 345--364.

\bibitem[Ful69]{fulton_on_hurwitz}
W.~Fulton, \emph{Hurwitz schemes and irreducibility of moduli of algebraic
  curves}, Ann. of Math. (2) \textbf{90} (1969), 542--575.

\bibitem[GHMS02]{GHS_converse}
T.~Graber, J.~Harris, B.~Mazur, and J.~Starr, \emph{Rational connectivity and
  sections of families over curves}, pre-print AG/0210225, 2002.

\bibitem[GHS03]{GHS}
T.~Graber, J.~Harris, and J.~Starr, \emph{Families of rationally connected
  varieties}, J. Amer. Math. Soc. \textbf{16} (2003), no.~1, 57--67
  (electronic).

\bibitem[Har77]{hartshorne}
R.~Hartshorne, \emph{Algebraic geometry}, Springer-Verlag, New York, 1977,
  Graduate Texts in Mathematics, No. 52.

\bibitem[HGS02]{GHS_Hurwitz}
J.~Harris, T.~Graber, and J.~Starr, \emph{A note on hurwitz schemes of covers
  of a positive genus curve}, pre-print AG/0205056, 2002.

\bibitem[HT04]{hassett_tschinkel}
B.~Hassett and Y.~Tschinkel, \emph{Weak approximation over function fields},
  pre-print AG/0411367, 2004.

\bibitem[KMM92a]{kmm3}
J.~Koll{\'a}r, Y.~Miyaoka, and S.~Mori, \emph{Rational connectedness and
  boundedness of {F}ano manifolds}, J. Differential Geom. \textbf{36} (1992),
  no.~3, 765--779.

\bibitem[KMM92b]{kmm2}
\bysame, \emph{Rationally connected varieties}, J. Algebraic Geom. \textbf{1}
  (1992), no.~3, 429--448.

\bibitem[Kol96]{kollar}
J.~Koll{\'a}r, \emph{Rational curves on algebraic varieties}, Ergebnisse der
  Mathematik und ihrer Grenzgebiete, vol.~32, Springer-Verlag, Berlin, 1996.

\bibitem[Kol01]{kollar_simplest_alg_var}
\bysame, \emph{Which are the simplest algebraic varieties?}, Bull. Amer. Math.
  Soc. (N.S.) \textbf{38} (2001), no.~4, 409--433 (electronic).

\bibitem[Kol04]{kollar_specialization}
\bysame, \emph{Specialization of zero cycles}, Publ. Res. Inst. Math. Sci.
  \textbf{40} (2004), no.~3, 689--708.

\end{thebibliography}

\end{document}